\documentclass[a4paper]{article}
\usepackage{amsfonts}
\usepackage{amsmath}
\usepackage{amssymb}
\begin{document}
\newtheorem{theorem}{Theorem}[section]
\newtheorem{lemma}[theorem]{Lemma}
\newtheorem{corollary}[theorem]{Corollary}
\newtheorem{conjecture}[theorem]{Conjecture}
\newtheorem{remark}[theorem]{Remark}
\newtheorem{definition}[theorem]{Definition}
\newtheorem{problem}[theorem]{Problem}
\newtheorem{example}[theorem]{Example}
\newtheorem{proposition}[theorem]{Proposition}
\title{{\bf Dynamical construction of K\"{a}hler-Einstein metrics}}
\date{November 24, 2006}
\author{Hajime TSUJI}
\maketitle
\begin{abstract}
\noindent In this article, we give a new construction of a K\"{a}hler-Einstein metric
on  a smooth projective variety with ample canonical bundle. 
This result can be generalized to the construction of a singular K\"{a}hler-Einstein metric on a smooth projective variety of general type which gives 
an AZD of the canonical bundle.

As a consequence, for a proper projective morphism $f : X \longrightarrow S$ 
(with connected fibers) such that a general fiber is of general type
and a positive integer $m$,
we construct a canonical singular hermitian metric $h_{E,m}$ on $f_{*}{\cal O}_{X}(mK_{X/S})$ with semipositive curvature in the sense of Nakano. 
\end{abstract}
\tableofcontents
\section{Introduction}

Let $X$ be a smooth projective $n$-fold with ample canonical bundle defined
over $\mathbb{C}$.  
Then by the celebrated solution of Calabi's conjecture (\cite{a,y}), there exists 
a unique K\"{a}hler-Einstein form $\omega_{E}$ such that 
\[
-\mbox{Ric}_{\omega_{E}} = \omega_{E}, 
\]
where $\mbox{Ric}_{\omega_{E}}$ denotes the Ricci form of the K\"{a}hler 
manifold $(X,\omega_{E})$. \vspace{3mm}\\

On the other hand for a complex manifolds with very ample $L^{2}$ canonical 
forms, there exists a standard K\"{a}hler form called the Bergman K\"{a}hler
form. 

Let us explain more precisely.  
Let $M$ be a complex manifold of dimension $n$ such that 
the space of $L^{2}$ canonical forms
\[
H^{0}_{(2)}(M,{\cal O}_{M}(K_{M})) := 
\{ \eta\in H^{0}(M,{\cal O}_{M}(K_{M})) \mid (\sqrt{-1})^{n^{2}}\int_{M}\eta\wedge\bar{\eta} < \infty \}
\]
gives a very ample linear system.
Then $M$ admits a Bergman kernel, 
\[
B(z,w) := \sum_{i} \sigma_{i}(z)\cdot\overline{\sigma_{i}(w)} \hspace{5mm}
(z,w\in M), 
\]
where $\{\sigma_{i}\}$ is a complete orthonormal basis of 
$H^{0}_{(2)}(M,{\cal O}_{M}(K_{M}))$ with respect to the inner product;
\[
(\eta,\eta^{\prime}):= (\sqrt{-1})^{n^{2}}\int_{M}\eta\wedge\bar{\eta}^{\prime}.
\] 
And
\[
\omega_{B}(z) := \sqrt{-1}\partial\bar{\partial}\log B(z,z) \hspace{5mm} (z\in M)
\]
is a K\"{a}hler form which is called the Bergman K\"{a}hler form on $M$.
The same construction applies for the case of the adjoint bundle of 
a (possibly singular) hermitian line bundle $(L,h)$ on $M$ (see Section 3). 
 
Both K\"{a}hler-Einstein metrics and Bergman metrics are determined uniquely
by the complex structures.   In this sense these metrics are canonical.
Hence it is natural to study the relation of these metrics. 
\vspace{3mm}\\

Recently S.K. Donaldson found a new construction of 
K\"{a}hler-Einstein metrics or more generally  K\"{a}hler metrics 
with constant scalar curvature. 
Actually he found a strong connection between the existence
of K\"{a}hler metrics with constant scalar curvature and the asymptotic 
stability of Hilbert points of projective embeddings (\cite{do}). 
In particular this implies the connection between the existence of 
K\"{a}hler-Einstein metrics and the asymptotic stability of Hillbert
points of projective embeddings (\cite{do}).
\vspace{3mm} \\
Let us explain (a part of) his results.
Let $X$ be a smooth projective variety and let $L$ be an ample line bundle
on $X$.
Then for every sufficiently large positive integer $m$, the linear system
$\mid\! mL\!\mid$ gives a projective embbedding
\[
\Phi_{m} : X \longrightarrow \mathbb{P}^{N_{m}}
\]
given by 
\[
\Phi_{m}(x) := [\sigma_{0}^{(m)}:\cdots :\sigma_{N_{m}}^{(m)}],
\]
where $\{\sigma_{0}^{(m)},\cdots ,\sigma_{N_{m}}^{(m)}\}$ 
is a basis of $H^{0}(X,{\cal O}_{X}(mL))$. 
Hence $\Phi_{m}$ depends on the choice of the basis.
Let $\omega_{FS}$ denote the Fubini-Study K\"{a}hler form on 
$\mathbb{P}^{N_{m}}$. 
If for some choice of $\{\sigma_{0}^{(m)},\cdots ,\sigma_{N_{m}}^{(m)}\}$ 
the equality
\[
\int_{X}\frac{\sigma_{i}^{(m)}\cdot\bar{\sigma}_{j}^{(m)}}{\sum_{i=0}^{N_{m}}\mid\sigma_{i}^{(m)}\mid^{2}} (\Phi_{m}^{*}\omega_{FS})^{n} = \delta_{ij} 
\]
holds for every $0\leqq i,j\leqq N_{m}$ (i.e.,$\{\sigma_{0}^{(m)},\cdots ,\sigma_{N_{m}}^{(m)}\}$ is orthonormal with respect to the $L^{2}$-inner product 
with respect to the hermitian metric $(\sum_{i=0}^{N_{m}}\mid\sigma_{i}^{(m)}\mid^{2})^{-2} $ on $mL$  and the volume form $(\Phi_{m}^{*}\omega_{FS})^{n}$), the K\"{a}hler form 
\[
\omega_{m}:= \frac{1}{m}\Phi_{m}^{*}\omega_{FS}
\]
is called {\bf balanced (or critical)}.  
The Hilbert point of $\Phi_{m}(X)$ is stable, if and only if 
there exists a choice of the basis $\{\sigma_{0}^{(m)},\cdots ,\sigma_{N_{m}}^{(m)}\}$ such that $\Phi_{m}$ is balanced (\cite{z}). 
Donaldson's theorem is stated as follows.

\begin{theorem}(\cite[p.482, Theorem 3]{do})\label{do}
Let $X$ be a smooth projective variety and let $L$ be an ample line bundle
 on $X$.  Suppose that $\mbox{Aut}(X,L)$ is discrete. 
If $X$ admits a K\"{a}hler form $\omega$ cohomologous to $2\pi c_{1}(L)$
with constant scalar curvature, 
then for every sufficiently large $m$, $\Phi_{m}(X)$ is stable (this property
is called that $(X,L)$ is asymptotically stable). 
And the limit of the balanced K\"{a}hler forms $\{ \omega_{m}\}$ 
exists in $C^{\infty}$-topology and the limit is a K\"{a}hler form with constant scalar curvature. $\square$
\end{theorem}
In short Theorem \ref{do} gives a construction of a K\"{a}hler form 
with constant scalar curvature as the limit of s sequence of balanced K\"{a}hler forms.  And Theorem \ref{do} is closely related to the asymptotic
expansion of Bergman kernels (\cite{c,ze}). \vspace{2mm}\\

In this article, we shall give a new   construction of K\"{a}hler-Einstein forms with negative Ricci curvature as a limit of Bergman K\"{a}hler forms. 
The purpose of this article is to relate K\"{a}hler-Einstein forms and 
Bergman K\"{a}hler forms in the case of projective manifolds with 
ample canonical bundle  or more generally projective manifolds of 
general type. 
\vspace{3mm} \\ 

Let $X$ be a smooth projective $n$-fold with ample canonical bundle. 
Let $m_{0}$ be a positive integer such that : 
\begin{enumerate}
\item $\mid mK_{X}\mid$ is very ample for every $m \geqq m_{0}$,
\item For every pseudoeffective singular hermitian line bundle $(L,h_{L})$(cf. Definition \ref{pe} below),  ${\cal O}_{X}(m_{0}K_{X} + L)\otimes{\cal I}(h_{L})$
is globally generated. 
\end{enumerate}
The existence of such $m_{0}$ follows from  Nadel's vanishing theorem
(\cite[p.561]{n}).

Let $h_{m_{0}}$ be a $C^{\infty}$ hermitian metric on $m_{0}K_{X}$ with 
strictly positive curvature.  Suppose that we have constructed 
$K_{m}$ and the $C^{\infty}$ hermitian metric $h_{m}$ on $mK_{X}$.  Then we define
\[
K_{m+1}:= K(X,(m+1)K_{X},h_{m})
\]
and 
\[
h_{m+1}:= 1/K_{m+1}, 
\]
where $K(X,(m+1)K_{X},h_{m})$ denotes (the diagonal part of) the Bergman kernel of $(m+1)K_{X}$
with respect to $h_{m}$ constructed as follows.
 
Let 
$\{\sigma_{0}^{(m+1)},\cdots ,\sigma_{N_{m+1}}^{(m+1)}\}$ be the 
complete orthonormal basis of \\ 
$H^{0}(X,{\cal O}_{X}((m+1)K_{X}))$  with respect to the inner product
\[
(\sigma,\tau ):= (\sqrt{-1})^{n_{2}}\int_{X}h_{m}\cdot\sigma\wedge\bar{\tau}\hspace{5mm}
(\sigma,\tau\in H^{0}(X,{\cal O}_{X}((m+1)K_{X}))).
\]
Then  for $x\in X$  we define
\begin{eqnarray*}
K_{m+1}(x)& = &K(X,(m+1)K_{X},h_{m})(x) \\
 &:= &\sum_{i=0}^{N_{m+1}}\mid\sigma_{i}^{(m+1)}\mid^{2}(x),
\end{eqnarray*}
where for a global section $\sigma$ of $(m+1)K_{X}$, 
$\mid \sigma\mid^{2}$ denotes the global section $\sigma\cdot\bar{\sigma}$
 of $(K_{X}\otimes\overline{K_{X}})^{\otimes (m+1)}$. 
We note that by the choice of $m_{0}$, $\mid\!(m+1)K_{X}\!\mid$ is very ample. 
Hence $h_{m+1} := 1/K_{m+1}$ is a $C^{\infty}$ hermitian metric on 
$(m+1)K_{X}$. 
Inductively we construct the sequences $\{ h_{m}\}_{m\geqq m_{0}}$
and $\{ K_{m}\}_{m > m_{0}}$.
This is the same construction originated by the author in \cite{tu3}. 

The following theorem is the main result in this article.

\begin{theorem}\label{main}
 Let $X$ be a smooth projective $n$-fold with ample canonical 
bundle.  Let $m_{0}$ and $\{ h_{m}\}_{m > m_{0}}$ be the sequence 
of hermitian metrics  as above.
Then 
\[
h_{\infty} := \liminf_{m\rightarrow\infty} \sqrt[m]{(m!)^{n}\cdot h_{m}}
\]
is a $C^{\infty}$ hermitian metric on $K_{X}$ such that 
\[
\omega_{\infty} := \sqrt{-1}\Theta_{h_{\infty}}
\]
is a K\"{a}hler form  on $X$ with 
\[
-\mbox{\em Ric}_{\omega_{\infty}} = \omega_{\infty}.
\]
$\square$
\end{theorem}
\begin{remark} The existence of the limit $h_{\infty}$ has already been proved 
in \cite{tu3}.  For the case of smooth projective varieties of non-general type  see \cite{tu3,tu5}. 
$\square$
\end{remark}

The construction of K\"{a}hler-Einstein form in Theorem \ref{main}
is more straightforward than the one in Theorem \ref{do}.
And Theorem \ref{main} seems to imply that the sequence of K\"{a}hler 
forms 
\[
\{\frac{\sqrt{-1}}{m}\Theta_{h_{m}}\}_{m\geqq m_{0}}
\]
induced by the projective morphisms $\Phi^{(m)} : X \longrightarrow 
\mathbb{P}^{N_{m}} (m > m_{0})$ defined by 
\[
\Phi^{(m)} (x) = [\sigma^{(m)}_{0}(x) :\cdots :\sigma_{N_{m}}^{(m)}(x)]
\hspace{5mm}(x \in X) 
\]
is asymptotically nearly balanced. 

We can generalize 
Theorem \ref{main} to the case of the maximal K\"{a}hler-Einstein current 
(cf. Definition \ref{maximal}) on a smooth projective variety of general type whose canonical bundle is not necessarily ample (Theorem \ref{main2}) without any essential change.
And this immediately implies the uniqueness of the maximal K\"{a}hler-Einstein 
currents  on smooth projective varieties of general type
(Theorem \ref{KE}).  
 
This enables us to deduce the logarithmic plurisubharmonicity of K\"{a}hler-Einstein volume forms
on a projective family (cf. Theorem \ref{positivity}) by using the 
recent result on variation of Bergman kernels (\cite{b1,b2,b3,tu5}).
And we are able to  study the degeneration of K\"{a}hler-Einstein 
currents on projective families of varieties of general type. 

\begin{theorem}\label{positivity}
Let $f : X \longrightarrow S$ be a  proper projective morphism with connected fibers
betweem smooth varieties. 
Let $S^{\circ}$ denote the maximal Zariski dense subset of $S$ such that 
$f$ is smooth over $X^{\circ}:= f^{-1}(S^{\circ})$. 
Suppose that a general fiber of $f$ is a smooth projective variety of 
general type. 
Let $\omega_{E/S}$ be the family of relative maximal K\"{a}hler-Einstein 
currents
on $X^{\circ}$ (cf. Definition \ref{maximal}).  Let $h^{\circ}_{E}$ be 
the singular hermitian metric on $K_{X/S}\!\mid\!X^{\circ}$ defined by
\[
h^{\circ}_{E}:= (\omega_{E/S}^{n})^{-1},
\]
where $n$ denotes the relative dimension of $f : X \longrightarrow S$. 
Then we have the followings :
\begin{enumerate}
\item $h_{E}^{\circ}$ extends to  a singular hermitian metric $h_{E}$ on $K_{X/S}$.
\item The curvature current $\Theta_{h_{E}}$ of $h_{E}$ is semipositive on $X$.
\end{enumerate}
$\square$ 
\end{theorem}
Theorem \ref{positivity} implies the following refinement of Kawamata's 
positivity theorem  (\cite[p.57, Theorem 1]{ka1}) for the direct image of 
 the relative pluricanonical bundle in the 
case that a general fiber is a variety of general type.

\begin{theorem}\label{nakano2}
Let $f : X \longrightarrow S$ be a proper projective morphism  with connected fibers betweem smooth varieties. 
Let $S^{\circ}$ denote the maximal Zariski dense subset of $S$ such that 
$f$ is smooth over $X^{\circ}:= f^{-1}(S^{\circ})$. 
Suppose that a general fiber of $f$ is a smooth projective variety of 
general type. 
Let $h_{E}$ be the singular hermitian metric on $K_{X/S}$ as in 
Theorem \ref{positivity}. 

$F_{m}:= f_{*}{\cal O}_{X}(mK_{X/S})$ is locally free on $S^{\circ}$ and 
$F_{m}\mid S^{\circ}$ carries the continuous hermitian metric 
$h_{E,m}$ defined by 
\[
h_{E,m}(\sigma ,\tau ):= (\sqrt{-1})^{n^{2}}\int_{X_{s}}h_{E}^{m-1}\cdot \sigma\wedge\bar{\tau} \hspace{5mm} (\sigma ,\tau\in H^{0}(X_{s},{\cal O}_{X_{s}}(mK_{X_{s}}))),
\]
where $n$ denotes the relative dimension of $f : X\longrightarrow S$. 
Then we have the followings.
\begin{enumerate}
\item The curvature $\Theta_{h_{E,m}}$ of $h_{E,m}$ is semipositive in the sense of Nakano.  
\item Let $x\in S- S^{\circ}$ be a point and let $\sigma$ be a local holomorphic
section of $F_{m}$ on a neighbourhood $U$ of $x$. 
Then $\sqrt{-1}\bar{\partial}\partial\log h_{E,m}(\sigma ,\sigma )$ extends to a closed positive current across 
$(S - S^{\circ})\cap U$.
\end{enumerate}
$\square$ 
\end{theorem}
\begin{remark} If $K_{X/S}$ is relatively ample, then $h_{E,m}$ is $C^{\infty}$
on $S^{\circ}$. $\square$
\end{remark}

\noindent Theorem \ref{positivity} has several applications. 
For example it  immediately gives canonical positive line bundles on the moduli space
of canonically polarized varieties with only canonical singularities. 
 Such applications will be discussed in the subsequent papers 
 because of the length.  

We should note that the convergence in Theorem \ref{main} is 
much weaker than in Theorem \ref{do}.  And Theorem \ref{main} does not 
say anything about K\"{a}hler forms with constant scalar curvature
at the moment.

\section{Preliminaries}

In this section, we shall review the basic terminologies used in this paper.  

\subsection{Singular hermitian metrics}\label{singh}
In this subection $L$ will denote a holomorphic line bundle on a complex manifold $M$. 
\begin{definition}\label{singhm}
A  {\bf singular hermitian metric} $h$ on $L$ is given by
\[
h = e^{-\varphi}\cdot h_{0},
\]
where $h_{0}$ is a $C^{\infty}$ hermitian metric on $L$ and 
$\varphi\in L^{1}_{loc}(M)$ is an arbitrary function on $M$.
We call $\varphi$ a  weight function of $h$. $\square$ 
\end{definition}
The curvature current $\Theta_{h}$ of the singular hermitian line
bundle $(L,h)$ is defined by
\[
\Theta_{h} := \Theta_{h_{0}} + \sqrt{-1}\partial\bar{\partial}\varphi ,
\]
where $\partial\bar{\partial}$ is taken in the sense of a current.
The $L^{2}$-sheaf ${\cal L}^{2}(L,h)$ of the singular hermitian
line bundle $(L,h)$ is defined by
\[
{\cal L}^{2}(L,h)(U) := \{ \sigma\in\Gamma (U,{\cal O}_{M}(L))\mid 
\, h(\sigma ,\sigma )\in L^{1}_{loc}(U)\} ,
\]
where $U$ runs over the  open subsets of $M$.
In this case there exists an ideal sheaf ${\cal I}(h)$ such that
\[
{\cal L}^{2}(L,h) = {\cal O}_{M}(L)\otimes {\cal I}(h)
\]
holds.  We call ${\cal I}(h)$ the {\bf multiplier ideal sheaf} of $(L,h)$.
If we write $h$ as 
\[
h = e^{-\varphi}\cdot h_{0},
\]
where $h_{0}$ is a $C^{\infty}$ hermitian metric on $L$ and 
$\varphi\in L^{1}_{loc}(M)$ is the weight function, we see that
\[
{\cal I}(h) = {\cal L}^{2}({\cal O}_{M},e^{-\varphi})
\]
holds.
For $\varphi\in L^{1}_{loc}(M)$ we define the multiplier ideal sheaf of $\varphi$ by 
\[
{\cal I}(\varphi ) := {\cal L}^{2}({\cal O}_{M},e^{-\varphi}).
\] 
\begin{example}
Let $\sigma\in \Gamma (X,{\cal O}_{X}(L))$ be the global section. 
Then 
\[
h := \frac{1}{\mid\sigma\mid^{2}} = \frac{h_{0}}{h_{0}(\sigma ,\sigma)}
\]
is a singular hermitian metric on $L$, 
where $h_{0}$ is an arbitrary $C^{\infty}$-hermitian metric on $L$
(the right hand side is obviously independent of $h_{0}$).
The curvature $\Theta_{h}$ is given by
\[
\Theta_{h} = 2\pi\sqrt{-1}(\sigma )
\]
where $(\sigma )$ denotes the current of integration over the 
divisor of $\sigma$. $\square$ 
\end{example}
\begin{definition}\label{pe}
$L$ is said to be {\bf pseudoeffective}, if there exists 
a singular hermitian metric $h$ on $L$ such that 
the curvature current 
$\Theta_{h}$ is a closed positive current.
Also a singular hermitian line bundle $(L,h)$ is said to be {\bf pseudoeffective}, 
if the curvature current $\Theta_{h}$ is a closed positive current. $\square$
\end{definition}

\subsection{Analytic Zariski decompositions}\label{AZD}

Let $L$ be a pseudoeffective line bundle on a compact complex manifold $X$.
To analyze the ring :
\[
R(X,L) = \oplus_{m=0}^{\infty}H^{0}(X,{\cal O}_{X}(mL))
\]
it is sometimes useful to  introduce the notion of analytic Zariski decompositions. 

\begin{definition}\label{defAZD}
Let $M$ be a compact complex manifold and let $L$ be a holomorphic line bundle
on $M$.  A singular hermitian metric $h$ on $L$ is said to be 
an analytic Zariski decomposition, if the followings hold.
\begin{enumerate}
\item $\Theta_{h}$ is a closed positive current,
\item for every $m\geq 0$, the natural inclusion
\[
H^{0}(M,{\cal O}_{M}(mL)\otimes{\cal I}(h^{m}))\rightarrow
H^{0}(M,{\cal O}_{M}(mL))
\]
is an isomorphim. $\square$
\end{enumerate}
\end{definition}
\begin{remark} If an AZD exists on a line bundle $L$ on a smooth projective
variety $M$, $L$ is pseudoeffective by the condition 1 above. $\square$
\end{remark}

It is known that for every pseudoeffective line bundle on a compact complex manifold, there exists an AZD on $L$ (cf. \cite{tu,tu2,d-p-s}). 

The advantage of the AZD is that we can handle pseudoeffective line bundle 
$L$ on a compact complex manifold $X$  
as a singular hermitian  line bundle with semipositive curvature current
as long as we consider the ring $R(X,L)$. 
 
\section{Variation of adjoint line bundles}

In this section we shall review the results in \cite{tu5}. 

\subsection{Theorems of Maitani-Yamaguchi and Berndtsson}

In 2004, Maitani and Yamaguchi proved the following theorem.

\begin{theorem}(\cite{m-y})\label{m-y} Let $\Omega$ be a pseudoconvex domain 
in $\mathbb{C}_{z}\times \mathbb{C}_{w}$ with $C^{1}$ boundary. 
Let $\Omega_{t} := \Omega \cap (\mathbb{C}_{z}\times \{ t\})$ and 
Let $K(z,t)$ be the Bergman kernel function of 
$\Omega_{t}$. 

Then $\log K(z,t)$ is a plurisubharmonic function on $\Omega$. $\square$
\end{theorem}
Recently  generalizing Theorem \ref{m-y}, B. Berndtsson proved the following
higher dimensional and twisted version of Theorem \ref{m-y}.
 
\begin{theorem}(\cite{b1})\label{b1}
Let $D$ be a pseudoconvex domain in $\mathbb{C}^{n}_{z}\times \mathbb{C}^{k}_{t}$.   And let $\phi$ be a plurisubharmonic function on $D$.
For $t\in \Delta$, we set  $D_{t} := \Omega \cap (\mathbb{C}^{n}\times \{ t\})$
and $\phi_{t} := \phi\mid D_{t}$. 
Let $K(z,t) (t\in \mathbb{C}^{k}_{t})$ be the Bergman kernel of the Hilbert space 
\[
A^{2}(D_{t},e^{-\phi_{t}}) := \{ f\in {\cal O}(\Omega_{t})\mid 
\int_{D_{t}}e^{-\phi_{t}}\mid f\mid^{2} < + \infty \} .
\]
Then $\log K(z,t)$ is a plurisubharmonic function on $D$. $\square$
\end{theorem}
As in mensioned in \cite{b2}, his proof also works for a pseudoconvex 
domain in a locally trivial family of manifolds which admits 
a Zariski dense Stein subdomain.  

Also he proved the following theorem. 

\begin{theorem}(\cite[Theorem 1.1]{b2})\label{b2}
Let us consider a domain $D = U \times\Omega$ and let $\phi$ be a 
plurisubharmonic function on $D$.
For simplicity we assume that $\phi$ is smooth up to the boundary and 
strictly plurisubharmonic in $D$. 
Then for each $t\in U$, $\phi_{t} := \phi (\cdot ,t)$ is plurisubharmonic 
on $\Omega$.  Let $A^{2}_{t}$ be the Bergman space of holomorphic functions
on $\Omega$ with norm
\[
\parallel f\parallel^{2} = \parallel f\parallel^{2}_{t}:=  \int_{\Omega} e^{-\phi_{t}}\mid f\mid^{2}. 
\]
The spaces $A^{2}_{t}$ are all equal as vector spaces but have norms 
that vary with $t$. 
Then ``infinite rank'' vector bundle $E$ over $U$ with fiber
$E_{t} = A^{2}_{t}$ is therefore trivial as a bundle but is equipped with 
a notrivial metric. 
Then $(E,\parallel \,\,\,\,\parallel_{t})$ is strictly positive 
in the sense of Nakano. $\square$ 
\end{theorem}
In Theorem \ref{b1} the assumption that $D$ is a pseudoconvex domain 
in the product space is rather strong.  And in Theorem \ref{b2}, Berndtsson
also assumed that $D$ is a product. 

\subsection{Variation of hermitian adjoint bundles}

Recently Berndtsson and I have independently generalized Theorems \ref{b1} and \ref{b2} to the case of projective families. 
The strategis of the proofs of \cite{b3} and \cite{tu5} are
completely
 different. 
Although \cite{tu5} depends heavily on \cite{b1,b2}, the result of \cite{tu5},
is substantially stronger than that of \cite{b3}. 
In fact in \cite{b3}, Berndtsson only consider the case of smooth fibrations,
but in \cite{tu5}, the method applies to the case of general projective 
fibrations which posseses singular fibers.

\begin{theorem}\label{v}(\cite{tu5})
Let $f : X \longrightarrow S$ be a projective family 
of projective varieties over a complex manifold $S$.
Let $S^{\circ}$  be the maximal nonempty Zariski open subset such that 
$f$ is smooth over $S^{\circ}$.

Let $(L, h)$ be a singular hermitian line bundle on $X$ such that 
$\Theta_{h}$ is semipositive on $X$. \\
Let $K_{s} := K(X_{s},K_{X}+ L\mid_{X_{s}},h\mid_{X_{s}})$ be the Bergman kernel 
of $K_{X_{s}}+ (L\mid X_{s})$ with respect to $h\mid X_{s}$ for $s\in S^{\circ}$. 
Then the singular hermitian metric $h_{B}$ of $K_{X/S} + L\mid f^{-1}(S^{\circ})$ defined by 
\[
h_{B}\mid X_{s}:= K_{s}^{-1} (s\in S^{\circ})
\]
has semipositive curvature on $f^{-1}(S^{\circ})$ and extends on $X$ as a singular hermitian metric on $K_{X/S} +L$ with semipositive curvature current. $\square$ 
\end{theorem}
\begin{theorem}\label{nakano}(\cite{tu5})
Let $f : X \longrightarrow S$ be a  projective family 
of over a complex manifold $S$
such that $X$ is smooth. 
Let $S^{\circ}$  be a nonempty Zariski open subset such that 
$f$ is smooth over $S^{\circ}$.
Let $(L,h)$ be a hermitian line bundle on $X$ such that 
$\Theta_{h}$ is semipositive on $X$.
We define the hermitian metric $h_{E}$ on  $E:= f_{*}{\cal O}_{X}(K_{X/S} +L)\mid S^{\circ}$
 by 
\[
h_{E}(\sigma ,\tau ) := (\sqrt{-1})^{n^{2}}\int_{X_{s}}h\cdot \sigma\wedge\bar{\tau} \hspace{5mm} (\sigma ,\tau \in H^{0}(X_{s},{\cal O}_{X_{s}}(K_{X_{s}}+L\mid_{X_{s}}))),
\] 
where $n$ denotes the relative dimension of $f : X \longrightarrow S$. 
Let $S_{0}$ be the maximal  Zariski open subset of $S$ such that 
$E\mid S_{0}$ is locally free. 
Then $(E,h_{E})\!\!\mid\!\!S_{0}$ is semipositive in the  sense of Nakano. 
Moreover if $\Theta_{h}$ is strictly positive, then 
$(E,h_{E})\!\!\mid\!S_{0}$ is strictly positive in the sense of Nakano.

Let $t\in S- (S^{\circ}\cap S_{0})$ be a point and let $\sigma$ be a local nonvanishing holomorphic
section of $E$ on a neighbourhood $U$ of $t$. 
Then $\sqrt{-1}\bar{\partial}\partial\log h_{E}(\sigma ,\sigma )$ extends as a closed positive current across 
$(S - (S^{\circ}\cap S_{0}))\cap U$. 
$\square$ 
\end{theorem}
{\em Proof of Theorems \ref{v} and \ref{nakano}}. \vspace{3mm} \\  
Let $f : X \longrightarrow S$ be a projective family. 
Since the statement is local we may assume that $S$ is the unit open ball 
$B$ with center $O$ in {$\mathbb{C}^{m}$. 
We may also and do assume that the family $f : X \longrightarrow B$ is 
a restriction of a projective family 
\[
\hat{f} : \hat{X} \longrightarrow B(O,2)
\]
over the open ball $B(O,2)$ of radius $2$ with center $O$.
Let 
\[
F : \hat{X} \times B(O,2) \longrightarrow B(O,2)\times B(O,2)
\]
the fiber space defined by 
\[
F(x,t) = (f(x),t). 
\]
Let $\varepsilon$ be a positive number less than $1$.
We set   
\[
T(\varepsilon ) = \{ (s,t) \in B(O,2)\times B(O,2) \mid 
t\in B, s\in B(t,\varepsilon )\}
\]
and  
\[
X(\varepsilon ) := F^{-1}(T(\varepsilon )).
\]
Let 
\[
f_{\varepsilon} : X(\varepsilon ) \longrightarrow B
\]
be the family  defined by 
\[
f_{\varepsilon}(x,t) = t.
\]
Since for $(x,t)\in X(\varepsilon )$, $x\in f^{-1}(B(t,\varepsilon ))$ holds,
we see that  
\[
X(\varepsilon )_{t} := f^{-1}(B(t,\varepsilon ))
\]
holds.  Hence we may consider $X(\varepsilon )_{t}$ as a family 
\[
\pi_{\varepsilon ,t}: X(\varepsilon ,t)\longrightarrow B(t,\varepsilon ).
\] 
We note that $T(\varepsilon )$ is a domain of holomorphy 
in $\mathbb{C}^{2m}$.  Hence $X(\varepsilon )$ is a pseudoconvex domain in $X\times B(O,2)$.  Since $X \times B(O,2)$ is a product manifold,  the proof of Theorem \ref{b1} works without any essential change in this case (cf. \cite{b1}).
Hence if we define $K_{\varepsilon}$ by 
\[
K_{\varepsilon}\mid X(\varepsilon )_{t}:= K(X(\varepsilon )_{t},K_{X}+L\mid X(\varepsilon )_{t},h_{L}\mid X(\varepsilon )_{t}), 
\]
then 
\begin{equation}\label{psh}
\sqrt{-1}\partial\bar{\partial}\log K_{\varepsilon}
 \geqq 0
\end{equation}
holds on $X(\varepsilon )$. 
We note that 
\begin{equation}\label{limit}
\lim_{\varepsilon\downarrow 0}\mbox{vol}(B(t,\varepsilon ))\cdot K(X(\varepsilon )_{t},K_{X}+L\mid\!X(\varepsilon )_{t},h_{L}\mid\!X(\varepsilon )_{t})
= K(X_{t},K_{X_{t}} + L\mid X_{t}, h\mid X_{t})
\end{equation}
holds.  
In fact, if we consider the family 
\[
\pi_{\varepsilon ,t} : X(\varepsilon ) \longrightarrow B(t,\varepsilon )
\] 
as a family over the unit open ball  $B$ in $\mathbb{C}^{m}$ with center $O$  by 
\[
t^{\prime}\mapsto \varepsilon^{-1}(t^{\prime} - t),
\]
the limit as $\varepsilon \downarrow 0$ is nothing but the 
trivial family $X_{t}\times B$.   
We note that  for a $L$-valued canonical form $\sigma$ on $f^{-1}(B(t,\varepsilon ))$,
\[
\int_{B(t,\varepsilon )}h_{E}(\sigma,\sigma)
\]
is nothing but the $L^{2}$-norm of the $L$-valued canonical form $\sigma$ with respect $h_{L}$ over $f^{-1}(B(t,\varepsilon ))$ by Fubini's theorem, 
where we abbrebiate the standard Lebesgue measure on $\mathbb{C}^{m}$. 
 Then the desired equality  follows from 
the $L^{2}$-extension theorm (\cite{o-t,o}) and the extremal property of the 
Bergman kernels. 

Combining (\ref{psh}) and (\ref{limit}), we complete the proof of  Theorem \ref{v}.

The proof of Theorem \ref{nakano}, is quite similar. 
First we note that 
\[
f_{\varepsilon} : X(\varepsilon ) \longrightarrow B
\]
is everywhere smooth.  

For the moment,  we shall assume that $E$ is locally free on $B(O,2)$.
Then there exists a  global generator $\{ \sigma_{1},\cdots ,\sigma_{r}\}$
of $E$ on $B(O,2)$, where $r = \mbox{rank}\, E$. 
Then we see that every $t\in B$, the fiber of the vector bundle 
$(f_{\varepsilon})_{*}{\cal O}_{X(\varepsilon )}(K_{X(\varepsilon )/B}+L)\otimes{\cal I}(h_{L})$ at $t$ is  canonically isomorphic to 
$\mathbb{C}^{r}\times {\cal O}(B(t,\varepsilon ))$ in terms of the frame
$\{ \sigma_{1},\cdots ,\sigma_{r}\}$.  And moreover the space
${\cal O}(B(t,\varepsilon ))$ is canonically isomorphic to ${\cal O}(B(O,\varepsilon ))$ by the parallel translation.  
In this case 
\[
E_{\varepsilon}: = (f_{\varepsilon,(2)})_{*}{\cal O}_{X(\varepsilon)}(K_{X(\varepsilon )/B}\otimes p_{1}^{*}L\otimes {\cal I}(p_{1}^{*}h))
\]
is a vector bundle of infinite rank on $B$, where 
\[
p_{1}: X(\varepsilon ) \longrightarrow X
\]
denotes  the first projection
\[
p_{1}(x,t) = x,\hspace{5mm}(x,t)\in  X(\varepsilon )
\]
and $(f_{\varepsilon,(2)})_{*}{\cal O}_{X(\varepsilon)}(K_{X(\varepsilon )/B}\otimes p_{1}^{*}L\otimes {\cal I}(p_{1}^{*}h))$ denotes the direct image 
of $L^{2}$-holomorphic sections. 
 
By the same proof as Theorem \ref{b2}, we see that the curvature current of 
$h_{E,\varepsilon}$ is well defined everywhere on $B$ and is semipositive
in the sense of Nakano.  
Letting $\varepsilon$ tend to $0$, the curvature $\Theta_{h_{E_{\varepsilon}}}$ converges to the curvature $\Theta_{h_{E}}$ operating 
on $E_{t}\otimes {\cal O}_{B,t}$ in the obvious manner
for every $t$ such that $f$ is smooth over $t$. 
Hence $\Theta_{h_{E}}$ is semipositive in the sense of Nakano. 

Let $t\in S- (S^{\circ}\cap S_{0})$ be a point and let $\sigma$ be a local 
nonvanishing holomorphic
section of $E$ on a neighbourhood $U$ of $t$. 
Since the assertion is local, for the proof we may and do assume $\sigma$ is defined on $B(O,2)$.
The existence of the extension of the current  $\sqrt{-1}\bar{\partial}\partial\log h_{E}(\sigma ,\sigma )$ can be verified as follows. 
  
Let us begin the following lemma which follows from \cite[p.27,Corollary 7.3]{b-t}.

\begin{lemma}(\cite[Corollary 7.3]{b-t})\label{ext}
Let $\{ T_{k}\}_{k=1}^{\infty}$ is a sequence of closed positive $(1,1)$ current on the unit open disk $\Delta$ in $\mathbb{C}$. 
Let $T_{k} = \sqrt{-1}\partial\bar{\partial}\varphi_{k}$. 
Suppose that $\varphi_{k}$ is $L^{\infty}_{loc}$ on $\Delta$ and $\{ \varphi_{k}\}$ converges to a plurisubharmonic function 
$\varphi $  on the punctured disk $\Delta^{*} = \Delta - \{ O\}$. 
Then $\{ T_{k}\}$ converges to a closed positive current on $\Delta$. $\square$
\end{lemma}

\noindent 
Let us consider the current 
\[
T_{k} := \sqrt{-1}\bar{\partial}_{s}\partial_{s}\log \{\frac{1}{\mbox{vol}(B(s,1/k))}\int_{B(s,,1/k)}h_{E}(\sigma,\sigma)\}
\]
on $B$ for $k\geqq 1$, where the integration is taken with respect to 
the standard Lebesgue measure on $\mathbb{C}^{m}$. 
We note that since 
\[
\frac{1}{\mbox{vol}(B(s,1/k))}\int_{B(s,1/k)}h_{E}(\sigma,\sigma)
\]
is plurisuperharmonic 
\[
\{\frac{1}{\mbox{vol}(B(s,1/k))}\int_{B(s,1/k)}h_{E}(\sigma,\sigma)\}_{k=1}^{\infty}
\]
is monotone increasing.
By Lemma \ref{ext}, we see that  $\sqrt{-1}\bar{\partial}\partial\log h_{E}(\sigma ,\sigma )$  is canonically extended 
across $t\in B$  where $E$ is locally free (here $m= \dim B$ may be bigger than $1$, but we may use slicing by curves to apply Lemma \ref{ext}). 

Next we shall consider the extension across the point where $E$ is not 
locally free. 
In this case we just need to slice $B$ by complex curves passing through the 
curve.   Since every torsion free coherent sheaf over a curve is always 
locally free, we see that $\sqrt{-1}\bar{\partial}\partial\log h_{E}(\sigma ,\sigma )$ is canonically extended across 
all the points on  $B$.  The extension of $h_{B}$ in Theorem \ref{v} is similar. 

This completes the proof of Theorems \ref{v} and \ref{nakano}. $\square$ 

\section{Proof of Theorem \ref{main}}
Let $X$ be a smooth projective $n$-fold with ample canonical bundle.
Let $m_{0}$ be a positive integer such that : 
\begin{enumerate}
\item $\mid mK_{X}\mid$ is very ample for every $m \geqq m_{0}$,
\item For every pseudoeffective singular hermitian line bundle
$(L,h_{L})$, \\ ${\cal O}_{X}(m_{0}K_{X} + L)\otimes{\cal I}(h_{L})$
is globally generated. 
\end{enumerate}
Let $h_{m_{0}}$ be a $C^{\infty}$ hermitian metric on $m_{0}K_{X}$
with strictly positive curvature. 
Let $\{ h_{m}\}_{m\geqq m_{0}}$ and $\{ K_{m}\}_{m > m_{0}}$ be the sequences of hermitian metrics and Bergman kernels 
constructed as in Section 1, i.e., $\{ h_{m}\}_{m\geqq m_{0}}$ and $\{ K_{m}\}_{m > m_{0}}$ are defined 
inductively by 
\[
K_{m+1} = K(X,K_{X}+mK_{X},h_{m})
\]
and 
\[
h_{m+1} = 1/K_{m+1}. 
\]
\vspace{5mm}\\
Let $\omega_{E}$ be the  K\"{a}hler-Einstein form on $X$ such that 
\[
-\mbox{Ric}_{\omega_{E}} = \omega_{E}. 
\]
Let $dV_{E} = (n!)^{-1}\omega_{E}^{n}$ be the volume form associated 
with $(X,\omega_{E})$. 
\begin{lemma}\label{lower}
\[
\limsup_{m\rightarrow\infty}\sqrt[m]{(m!)^{-n}K_{m}}
\geqq (2\pi)^{-n}dV_{E}
\]
holds on $X$. $\square$ 
\end{lemma}
{\em Proof}. 
Let us consider the hermitian line bundle $(K_{X},dV_{E})$ on $X$.
Let $p\in X$ be a point.  Then by the K\"{a}hler-Einstein condition, there exists a holomorphic normal coordinate
$(U,z_{1},\cdots ,z_{n})$ such that 
\begin{equation}\label{(1)}
dV_{E}^{-1}=  \{ \prod_{i=1}^{n}(1 -\mid z_{i}\mid^{2}) + O(\parallel z\parallel^{3})\}\cdot 
\mid dz_{1}\wedge\cdots\wedge dz_{n}\mid^{-2}
\end{equation} 
holds.
Suppose that
\[
C_{m-1}\cdot dV_{E}^{m-1}\leqq K_{m-1}
\]
holds on $X$ for some positive constant $C_{m-1}$. 
We note that
\begin{equation}\label{extremal}
K_{m}(x) = \sup \{ \mid\sigma\mid^{2}(x) ;  
\sigma \in H^{0}(X,{\cal O}_{X}(mK_{X})), (\sqrt{-1})^{n^{2}}\!\int_{X}h_{m-1}\cdot\sigma\wedge
\bar{\sigma} = 1 \}
\end{equation} 
holds for every $x\in X$, by the extremal property of the Bergman kernel
(This is well known. See for example, \cite[p.46, Proposition 1.3.16]{kr}). 
We note that  for the open unit disk $\Delta = \{ t\in \mathbb{C}\mid \,\,\mid t\mid < 1\}$, 
\begin{equation}\label{disk}
\int_{\Delta}(1 - \mid t\mid^{2})^{m}dt\wedge d\bar{t}
= \frac{2\pi}{m+1} 
\end{equation}
holds.  Then  by H\"{o}rmander's $L^{2}$-estimate of $\bar{\partial}$-operator, 
we see that there exists a positive constant $\lambda_{m}$ such that 
\begin{equation}\label{induction}
(\lambda_{m}\cdot(2\pi)^{-n}\cdot m^{n})\cdot C_{m-1}\cdot  dV_{E}^{m}\leqq K_{m}
\end{equation}
with 
\[
\lambda_{m} \geqq 1 - \frac{C}{\sqrt{m}},
\]
where $C$ is a positive constant independent of $m$. 
\vspace{3mm} \\
In fact this can be verified as follows. 
Let $x\in X$ be a point on $X$ and let $(U,z_{1},\cdots ,z_{n})$
be the normal coordinate as above.  We may assume that 
$U$ is biholomorphic to the polydisk $\Delta^{n}(r)$ of radius $r$ with center $O$
in $\mathbb{C}^{n}$ for some $r$ via $(z_{1},\cdots ,z_{n})$. 

Taking $r$ sufficiently small we may assume that 
there exists a  $C^{\infty}$ function  $\rho$ on $X$ such that 
\begin{enumerate}
\item $\rho$ is identically $1$ on $\Delta^{n}(r/3)$.
\item $0\leqq \rho \leqq 1$. 
\item $\mbox{Supp}\,\rho \subset\subset U$.
\item $\mid d\rho\mid < 3/r$, where $\mid\,\,\,\,\,\mid$ denotes 
the pointwise norm with respect to $\omega_{E}$. 
\end{enumerate}
We note that by  the equation (\ref{(1)}), the mass of 
$\rho\cdot (dz_{1}\wedge\cdots \wedge dz_{n})^{\otimes m}$ concentrates 
around the origin as $m$ tends to infinity. 
Hence  by (\ref{disk}) we see that the $L^{2}$-norm 
\[
\parallel\rho\cdot (dz_{1}\wedge\cdots \wedge dz_{n})^{\otimes m}\parallel
\]
 of $\rho\cdot (dz_{1}\wedge\cdots \wedge dz_{n})^{\otimes m}$
with respect to $(dV_{E})^{-\otimes m}$ and $\omega_{E}$ is asymptotically
\begin{equation}\label{(2)}
\parallel\rho\cdot (dz_{1}\wedge\cdots \wedge dz_{n})^{\otimes m}\parallel^{2}
\sim (\frac{2\pi}{m})^{n}
\end{equation}
as $m$ tends to infinity , where $\sim$ means that  the ratio of the both sides 
converges to $1$. 
We set 
\[
\phi := n\rho\log \sum_{i=1}^{n}\mid z_{i}\mid^{2}. 
\]
We may and do assume that $m$ is sufficiently large so that
\[
m\cdot\omega_{E} +\sqrt{-1}\partial\bar{\partial}\phi > 0
\]
holds on $X$. 

By (\ref{(2)}), the $L^{2}$-norm 
\[
\parallel\bar{\partial}(\rho\cdot (dz_{1}\wedge\cdots \wedge dz_{n})^{\otimes m})\parallel_{\phi}
\]
of 
$\bar{\partial}(\rho\cdot (dz_{1}\wedge\cdots \wedge dz_{n})^{\otimes m})$ 
with respect to $e^{-\phi}\cdot(dV_{E})^{-\otimes m}$ and $\omega_{E}$ 
satisfies the inequality 
\begin{equation}\label{(3)}
\parallel\bar{\partial}(\rho\cdot (dz_{1}\wedge\cdots \wedge dz_{n})^{\otimes m})\parallel_{\phi}^{2}\leqq C_{0}\cdot (\frac{3}{r})^{2n+2}(\frac{2\pi}{m})^{n}
\end{equation}
for every $m$, where $C_{0}$ is a positive constant independent of $m$. 

By  H\"{o}rmander's $L^{2}$-estimate applied to the adjoint line bundle of \\ the hermitian line bundle 
$((m-1)K_{X},e^{-\phi}\cdot dV_{E}^{-(m-1)})$, we see that for every 
sufficiently large $m$, there exists a $C^{\infty}$ solution of 
the equation ; 
\[
\bar{\partial}u = \bar{\partial}(\rho\cdot (dz_{1}\wedge\cdots \wedge dz_{n})^{\otimes m})
\]
such that 
\[
u(x) = 0
\]
and  
\[
\parallel u\parallel_{\phi}^{2} \leqq \frac{2}{m}\parallel\bar{\partial}(\rho\cdot (dz_{1}\wedge\cdots \wedge dz_{n})^{\otimes m})\parallel_{\phi}^{2}
\]
hold, where $\parallel\,\,\,\,\parallel_{\phi}$'s  denote the $L^{2}$ norms
with respect to $e^{-\phi}\cdot dV_{E}^{-(m-1)}$ and $\omega_{E}$
respectively. 
Then $\rho\cdot (dz_{1}\wedge\cdots \wedge dz_{n})^{\otimes m}- u$
is a holomorphic section of $mK_{X}$ such that 
\[
(\rho\cdot (dz_{1}\wedge\cdots \wedge dz_{n})^{\otimes m}- u)(x)
= (dz_{1}\wedge\cdots \wedge dz_{n})^{m}
\]
and
\[
\parallel\rho\cdot (dz_{1}\wedge\cdots \wedge dz_{n})^{\otimes m}- u\parallel^{2}
 \leqq (1 + C_{0}\cdot (\frac{3}{r})^{2n+2}\sqrt{\frac{2}{m}})(\frac{2\pi}{m})^{n}.
\]
 
Hence by induction on $m$, using (\ref{extremal}) and 
(\ref{induction}), we see that there exist  positive constants $C$ and $C^{\prime}$ such that for every $m >  m_{0}$
\[
K_{m} \geqq C ^{\prime}(\prod_{k=m_{0}}^{m}(1 -\frac{C}{\sqrt{k}}))\cdot (m!)^{n}\cdot (2\pi )^{-mn}\cdot dV_{E}^{m}
\]
holds on $X$.   This implies that  
\[
\limsup_{m\rightarrow\infty}\sqrt[m]{(m!)^{-n}K_{m}}
\geqq (2\pi)^{-n}dV_{E}
\]
holds on $X$.
$\square$

\begin{lemma}\label{holder}
\[
\int_{X}\sqrt[m]{K_{m}} \leqq (\prod_{k=m_{0}}^{m}(N_{k}+1))^{\frac{1}{m}}\cdot (\int_{X}\sqrt[m_{0}]{K_{m_{0}}})^{\frac{m_{0}}{m}}
\]
holds, where $N_{k} := \dim \mid kK_{X}\mid =\dim H^{0}(X,{\cal O}_{X}(kK_{X})) -1$. $\square$
\end{lemma}
{\em Proof}.
By H\"{o}lder's ineqality we have 
\begin{eqnarray*}
\int_{X}\sqrt[m]{K_{m}} & = & \int_{X}\frac{K_{m}^{\frac{1}{m}}}
{K_{m-1}^{\frac{1}{m-1}}}\cdot K_{m-1}^{\frac{1}{m-1}}
\\
& \leqq  & (\int_{X}\frac{K_{m}} 
{K_{m-1}^{\frac{m}{m-1}}}\cdot  K_{m-1}^{\frac{1}{m-1}})^{\frac{1}{m}}
\cdot (\int_{X}K_{m-1}^{\frac{1}{m-1}})^{\frac{m-1}{m}} \\
& = & (\int_{X}\frac{K_{m}} 
{K_{m-1}})^{\frac{1}{m}}
\cdot (\int_{X}K_{m-1}^{\frac{1}{m-1}})^{\frac{m-1}{m}} \\
& = & (N_{m}+1)^{\frac{1}{m}}
\cdot (\int_{X}K_{m-1}^{\frac{1}{m-1}})^{\frac{m-1}{m}} \\
\end{eqnarray*}
Then continuing this process, by using 
\[
\int_{X}K_{m-1}^{\frac{1}{m-1}} \leqq (N_{m-1}+1)^{\frac{1}{m-1}}\cdot
(\int_{X}K_{m-2}^{\frac{1}{m-2}})^{\frac{m-2}{m-1}},
\]
we have that 
\[
\int_{X}(K_{m})^{\frac{1}{m}} \leqq \{(N_{m}+1)\cdot (N_{m-1}+1)\}^{\frac{1}{m}}
\cdot (\int_{X}(K_{m-2})^{\frac{1}{m-2}})^{\frac{m-2}{m}}
\]
holds.
Continueing this process we obtain the lemma. $\square$ \vspace{3mm} \\ 

Using Lemma \ref{holder}, we obtain the following lemma. 
\begin{lemma}\label{upper}
\[
\limsup_{m\rightarrow\infty}\frac{1}{(m!)^{\frac{n}{m}}}\int_{X}(K_{m})^{\frac{1}{m}} \leqq \frac{K_{X}^{n}}{n!}
\]
holds. $\square$
\end{lemma}
{\em Proof}.  
By the Kodaira vanishing theorem, 
\[
H^{q}(X,{\cal O}_{X}(mK_{X})) = 0
\]
holds for every $m\geqq 2$ and $q \geqq 1$. 
Then by  the Hirzebruch Riemann-Roch theorem, we have that 
\[
N_{m}+1 = \frac{K_{X}^{n}}{n!}m^{n} + O(m^{n-1})
\]
holds. Then by Lemma \ref{holder}, we see that 
\[
\limsup_{m\rightarrow\infty}\frac{1}{(m!)^{\frac{n}{m}}}\int_{X}(K_{m})^{\frac{1}{m}} \leqq \frac{K_{X}^{n}}{n!}
\]
holds. $\square$ \vspace{5mm}\\ 

Combining Lemmas \ref{lower} and \ref{holder}, we have the equality,
\[
\limsup_{m\rightarrow\infty}\frac{1}{(m!)^{\frac{n}{m}}}\sqrt[m]{K_{m}}
= (2\pi)^{-n}dV_{E},
\]
since
\[
\int_{X}dV_{E} = \frac{1}{n!}\int_{X}\omega_{E}^{n}
= \frac{(2\pi )^{n}K_{X}^{n}}{n!}
\]
hold by the K\"{a}hler-Einstein condition. 
This completes the proof of Theorem \ref{main}. $\square$

\section{Dynamical construction of K\"{a}hler-Einstein currents}
In this section we shall generalize Theorem \ref{main} to the case of general 
smooth projective varieties of general type.

\subsection{Existence of K\"{a}hler-Einstein currents}\label{exist}

In this subsection, we shall review the existence of a K\"{a}hler-Einstein 
current on a smooth projective variety of general type. 
Allowing singularities, there are infinitely many choice of K\"{a}hler-Einstein metrics.  But we focus on the metrics with minimal singularities.

\begin{theorem}\label{KE}
Let $X$ be a smooth projectie $n$-fold of general type. 
Then there exists a closed positive current $\omega_{E}$ on $X$ 
such that 
\begin{enumerate}
\item $\omega_{E}$ is $C^{\infty}$ on a nonempty Zariski open subset 
$U$ of $X$.
\item $\omega_{E}$ is a K\"{a}hler-Einstein metric on $U$ with 
\[
\omega_{E}=-\mbox{\em Ric}_{\omega_{E}}
\]
on $U$. 
\item The singular hermitian metric $(\omega_{E}^{n})^{-1}$ is an 
AZD of $K_{X}$.
\end{enumerate}
$\square$
\end{theorem}
\begin{definition}\label{maximal}
Let $X$ be a smooth projectie $n$-fold of general type and 
let $\omega_{E}$ be the closed positive current on $X$ as in 
Theorem \ref{KE}.  
We call $\omega_{E}$ the {\bf maximal K\"{a}hler-Einstein current} on $X$. 
$\square$
\end{definition}
\begin{remark}
Later we will see that the maximal K\"{a}hler-Einstein current 
is unique. $\square$ 
\end{remark}
{\em Proof of Theorem \ref{KE}}. The proof is more or less parallel to that of \cite[Theorem 5.6, p.430]{s}.  
Let $m_{0}$ be a sufficiently large positive integer such that 
$\mid\!m!K_{X}\!\mid$ gives a birational embedding of $X$ for every 
$m \geqq m_{0}$. 
Let $\pi_{m} : X_{m}\longrightarrow X$ be the resolution of $\mbox{Bs}\mid m!K_{X}\mid$ such that for every $m > m_{0}$
\[
\pi_{m} : X_{m}\longrightarrow X
\]
factors through $\pi_{m-1} : X_{m-1}\longrightarrow X$. 
Let 
\[
\mu_{m} : X_{m} \longrightarrow X_{m-1}
\]
be the natural morphism.
Let 
\[
\pi_{m}^{*}\mid\! m!K_{X}\!\mid = \mid\! P_{m}\!\mid + F_{m}
\]
be the decomposition of $\pi_{m}^{*}\mid\! m!K_{X}\!\mid$ into the 
free part $\mid P_{m}\mid$ and the fixed component $F_{m}$. 
Let $V$ be an analytic subset of $X$ defined by 
\[
V:= \{ x\in X\mid \mbox{$\Phi_{\mid\! m!K_{X}\!\mid}$ is not an embedding 
on an neighbourhood of $x$ for some $m \geqq m_{0}$}\}.
\]
There exists an effective $\mathbb{Q}$-divisor 
$E_{m}$ on $X_{m}$ respectively  such that 
\begin{enumerate}
\item $P_{m}- E_{m}$ is ample on $X_{m}$. 
\item $\mbox{Supp}\,E_{m}$ is contained in $\pi_{m}^{-1}(V)$. 
\item $((m+1)!)^{-1}(P_{m+1}- E_{m+1})-\mu_{m+1}^{*}(m!)^{-1}(P_{m}- E_{m}))$
is effective. 
\end{enumerate}
hold for every $m\geqq m_{0}$. 
After taking such a sequence $\{ E_{m}\}$, we replace $\{ E_{m}\}$ by 
$\{ 2^{-m}E_{m}\}$. 
Then it has the same properties as above. 
And we shall denote $\{ 2^{-m}E_{m}\}$ again by $\{ E_{m}\}$. 

Then by \cite[Theorem 5.6]{s}, there exists a closed positive current
$\omega_{m}$ on $X_{m}$ such that 
\begin{enumerate}
\item $-\mbox{Ric}_{\omega_{m}} = \omega_{m}$  holds on $X_{m}- \mbox{Supp}\, E_{m}$,
\item The absolutely continuous part of $\omega_{m}$ represents
$2\pi (m!)^{-1}(P_{m}-E_{m})$
\item $(\pi_{m})_{*}\omega_{m}$ represents the class $2\pi c_{1}(K_{X})$,
\end{enumerate}
Let us consider $\{ \omega_{m}^{n}\}$ as a sequence of volume forms 
on $X - V$.  And we shall identify $(\pi_{m})_{*}\omega_{m}$ and 
$\omega_{m}$ on $X - V$. 
Then by the maximum principle we see that 
\[
\omega_{m}^{n} \leqq \omega_{m+1}^{n}
\]
holds on $X - V$, by using the Einstein condition. 

Let 
\[
\mu : Y \longrightarrow X
\]
be a modification such that $\mu^{-1}(V)$ is a divisor with normal crossings.
Let $H$ be a sufficiently ample divisor on $Y$ such that 
\[
D := \mu^{-1}(V) + H
\]
is a divisor with normal crossings and $K_{Y} + D$ is ample. 
By \cite{kob}, there exists a complete 
K\"{a}hler-Einstein form $\omega_{D}$ on $Y -D$. 
Let us consider $\omega_{D}$ as a complete K\"{a}hler form on $X - \mu(D)$.
On the other hand by Yau's Schwarz lemma (\cite{y2}), we see that 
\[
\omega_{m}^{n} \leqq \omega_{D}^{n}
\]
holds on $X - \mu(D)$ for every $m$. 
Hence by moving $D$, 
\[
\lim_{m\rightarrow\infty}\omega_{m}^{n}
\]
exists on $X - V$. 

Now we shall consider the uniform $C^{2}$-estimate on every compact 
subset of $X - V$.  
Let $F$ be an effective $\mathbb{Q}$-divisor on $X$ such that 
\begin{enumerate}
\item $K_{X} - F$ is ample.
\item $\pi_{m}^{*}F - E_{m}$ is effective and  
$\mbox{Supp}\, (\pi_{m}^{*}F - E_{m})$ contains $\mbox{Supp}\, E_{m}$. 
\end{enumerate}
The existence of such a divisor $F$ follows from the proof of Kodaira's lemma.
Let $H$ be a smooth very ample divisor on $X$. 
Considering the exact sequence
\[
0\rightarrow H^{0}(X,{\cal O}_{X}(\ell K_{X} - H)\otimes {\cal I}_{V})
\rightarrow H^{0}(X,{\cal O}_{X}(\ell K_{X})) 
\rightarrow H^{0}(X,{\cal O}_{X}(\ell K_{X})\otimes{\cal O}_{X}/{\cal I}_{H}\cdot{\cal I}_{V})
\]
for $\ell >> 1$, we may find an effective member 
$F^{\prime}\in \mid\!H^{0}(X,{\cal O}_{X}(\ell K_{X} - H))\!\!\mid$.
Then $K_{X}- \ell^{-1}F^{\prime}$ is ample.  
Then $F$ can be taken as $\ell^{-1}F^{\prime}$. 
By this argument we see that $\cap_{F}\mbox{Supp}\, F = V$
holds, where $F$ runs all such $F$'s. 
Let 
\[
F = \sum a_{i}F_{i}
\]
be the irreducible decomposition of $F$.
Let $\sigma_{i}$ be a global section of ${\cal O}_{X}(F_{i})$
with divisor $F_{i}$ respectively. 
Let $\Omega$ be a $C^{\infty}$-volume form on $X$. 
Then since $K_{X}-F$ is ample, there exist hermitian metrics
$\{ h_{i}\}$ of $\{{\cal O}_{X}(F_{i})\}$ respectively such that 
\[
\omega_{F}:= -\mbox{Ric}\,\Omega -\sqrt{-1}\sum_{i}a_{i}\Theta_{h_{i}}
\]
is a K\"{a}hler form on $X$. 
We note that $\pi_{m}^{*}F - E_{m}$ is effective. 
Let $u_{m}$ be a $C^{\infty}$-function on $X - D$ such that 
\[
\omega_{m} = \omega_{F} + \sqrt{-1}\partial\bar{\partial}u_{m}
\]
and 
\[
\log \frac{(\omega_{F} + \sqrt{-1}\partial\bar{\partial}u_{m})^{n}}{\Omega\cdot\prod_{i}\parallel\sigma_{i}\parallel^{2a_{i}}}= u_{m}
\]
hold. 
We note that $u_{m}$ is identically $+\infty$ on  $F$ by the choice of 
$F$.  
Hence there exists a point $p_{0}\in X - F$, where $u_{m}$ takes its minimum. 
Then 
\[
\sqrt{-1}\partial\bar{\partial}\log \frac{(\omega_{F} + \sqrt{-1}\partial\bar{\partial}u_{m})^{n}}{\Omega\cdot\prod_{i}\parallel\sigma_{i}\parallel^{2a_{i}}}(p_{0}) \geqq 0
\]
holds. 
Hence  
\[
\omega_{m}(p_{0}) - \omega_{F}(p_{0}) \geqq 0 
\]
holds. This implies that 
\[
u_{m}(p_{0}) \geqq   \log \frac{\omega_{F}^{n}}{\Omega\cdot\prod_{i}\parallel\sigma_{i}\parallel^{2a_{i}}}(p_{0})
\]
holds. Hence we see that 
\[
u_{m}(x) \geqq \log \frac{\omega_{F}^{n}}{\Omega\cdot\prod_{i}\parallel\sigma_{i}\parallel^{2a_{i}}}(p_{0})
\]
holds for every $x\in X - V$.

Similarly if $u_{m} + \sum_{i}2a_{i}\log \parallel\sigma_{i}\parallel$ takes
its  maximum at a point $p_{0}^{\prime}$ on $X - V$,
since
\[
\log \frac{\omega_{m}^{n}}{\Omega} = u_{m} + \sum_{i}2a_{i}\log \parallel\sigma_{i}\parallel
\]
holds. 
We note that since $\pi_{m}^{*}F - E_{m}$ is effective, 
\[
u_{m} + \sum_{i}2a_{i}\log \parallel\sigma_{i}\parallel = -\infty 
\]
holds on $\mbox{Supp}\,F$.  This means that such a point $p_{0}^{\prime}$ 
certainly exists.  At $p_{0}^{\prime}$ we have that 
\[
\sqrt{-1}\partial\bar{\partial}\log \frac{\omega_{m}^{n}}{\Omega}(p_{0}^{\prime}) \leqq 0
\]
holds.
Hence noting $\omega_{m}$ is K\"{a}hler-Einstein, we see that 
\[
\omega_{m}(p_{0}^{\prime}) \leqq  (-\mbox{Ric}\,\Omega )(p_{0}^{\prime})
\] 
holds
and  
\[
u_{m} + \sum_{i}2a_{i}\log \parallel\sigma_{i}\parallel 
\leqq \log \frac{(-\mbox{Ric}\,\Omega )^{n}}{\Omega}(p_{0}^{\prime})
\]
holds on $X$. 

By the above consideration we have the following lemma. 
\begin{lemma}\label{c0}There exists a positive 
constant $C_{0}$ independent of $m$ such that 
\[
-C_{0} \leqq u_{m}  \leqq C_{0} - \sum_{i}2a_{i}\log \parallel\sigma_{i}\parallel
\]
hold on $X-V$. 
$\square$
\end{lemma}
\begin{lemma}\label{c2}(\cite[p. 127, Lemma 2.2]{tu0}))
We set 
\[
f := \log \frac{\omega_{F}^{n}}{\Omega}.
\]
Let $C$ be a positive number such that 
\[
C + \inf_{i\neq\ell}R_{i\bar{i}\ell\bar{\ell}} > 1
\]
holds on $X$, where $R_{i\bar{i}\ell\bar{\ell}}$ denotes the 
bisectional curvature. 

Then 
\[
e^{Cu_{m}}\Delta_{m}(e^{-Cu_{m}}(n + \Delta_{F}\,u_{m}))
\geqq (n + \Delta_{F}\,u_{m}) 
\]
\[
 +  \Delta_{F}(f+ \sum_{i}2a_{i}\log\parallel\sigma_{i}\parallel ) +(n+n^{2}\inf_{i\neq\ell}R_{i\bar{i}\ell\bar{\ell}} ) 
\]
\[
+ C\cdot n(n+\Delta_{F}\,u_{m}) + 
(n+\Delta_{F}\,u_{m})^{\frac{n}{n-1}}\exp(-\frac{1}{n-1}u_{m} + f)
\]
holds, where $\Delta_{F}$ denotes the Laplacian with respect to $\omega_{F}$
\mbox{\em (}i.e., $\Delta_{F}= \mbox{\em trace}_{\omega_{F}}\sqrt{-1}\partial\bar{\partial}$\mbox{\em )} and 
$\Delta_{m}$ denotes the Laplacian with respect to $\omega_{m}$. 

$\square$
\end{lemma}
Let $x_{0}$ be the point where $e^{-Cu_{m}}(n + \Delta u_{m})$ takes its maximum.
Then 
\[
0\leqq n + \Delta_{F}\,u_{m}(x_{0})\leqq C_{2}
\]
holds.
\[
0\leqq n + \Delta_{F}\,u_{m} \leqq  \exp (C(u_{m}- u_{m}(x_{0}))\cdot C_{2}
\]
By Lemma \ref{c0}, there exists a positive constant $C_{3}$ such that 
\[
n + \Delta_{F}\,u_{m} \leqq C_{3}(\prod_{i}\parallel\sigma_{i}\parallel^{2a_{i}})^{-C}
\]
holds on $X - V$.

Applying the general theory of fully nonlinear elliptic equations (\cite{tr}), 
moving $F$, 
we get a uniform higher order estimate of $u_{m}$ on every compact 
subset of $X - V$.  Letting $m$ tend to infinity, we see that 
by the monotonicity of $\{\omega_{m}^{n}\}$ 
\[
\omega_{E}:= \lim_{m\rightarrow\infty}\omega_{m}
\]
exists in $C^{\infty}$-topology on every compact subset of $X- V$.
Then it is clear that $\omega_{E}$ is K\"{a}hler-Einstein.
By the construction of $\{ E_{m}\}$ and the monotonicity of 
$\{ \omega_{m}^{n}\}$, we see that $(\omega_{E}^{n})^{-1}$ is an AZD of 
$K_{X}$. 
This completes the proof of Theorem \ref{KE}. $\square$
\subsection{A generalization of Theorem \ref{main}}\label{m2}

Let $X$ be a smooth projective $n$-fold of general type whose canonical
bundle is not necessarily ample. 

In this case we may also define the dynamical system of Bergman kernels
as in the case that $X$ has ample canonical bundle. 
In fact the construction of the dynamical system of Bergman kernel is 
parallel except the following differences.
\begin{enumerate}
\item The starting line bundle is not a multiple of $K_{X}$, but a 
sufficiently ample line bundle.
\item The hermitian metrics $\{ h_{m}\}$ are singular.
\end{enumerate}
Let us explain in detail. 
Let $A$ be a very ample line bundle on $X$  such that 
for every pseudoeffective singular hermitian line bundle $(L,h_{L})$, \\ ${\cal O}_{X}(m_{0}K_{X} + L)\otimes{\cal I}(h_{L})$
is globally generated. 

The existence of such $A$ follows from Nadel's vanishing theorem (\cite[p.561]{n}).

Let $h_{0}$ be a $C^{\infty}$-hermitian metric on $A$ with 
strictly positive curvature.  
Let $\{\sigma_{0}^{(1)},\cdots ,\sigma_{N_{1}}^{(1)}\}$ 
be a complete orthonormal basis of 
$H^{0}(X,{\cal O}_{X}(A+K_{X}))$ with respect to the inner product
\[
(\sigma,\tau ):= (\sqrt{-1})^{n_{2}}\int_{X}h_{0}\cdot\sigma\wedge\bar{\tau}
\hspace{5mm}(\sigma,\tau\in H^{0}(X,{\cal O}_{X}(A+K_{X}))),
\]
where we have considered $\sigma$ and $\tau$ as  $A$-valued 
$(n,0)$ forms.
We set 
\[
K_{1} = \sum_{i=0}^{N_{1}}\mid\!\sigma_{i}^{(1)}\!\!\mid^{2}
\]
and 
\[
h_{1}:= 1/K_{1}.
\]
It is clear that $K_{1}$ is independent of the choice of 
the orthonormal basis $\{\sigma_{0}^{(1)},\cdots ,\sigma_{N_{1}}^{(1)}\}$.
Suppose that $h_{m}$ is defined for some $m\geqq m_{0}+1$.
Then we define $h_{m+1}$ as follows. 
Let 
$\{\sigma_{0}^{(m+1)},\cdots ,\sigma_{N_{m+1}}^{(m+1)}\}$ be a  
complete orthonormal basis of 
$H^{0}(X,{\cal O}_{X}(A+(m+1)K_{X}))$  with respect to the inner product
\[
(\sigma,\tau ):= (\sqrt{-1})^{n_{2}}\int_{X}h_{m}\cdot\sigma\wedge\bar{\tau}
\hspace{5mm} 
(\sigma,\tau\in H^{0}(X,{\cal O}_{X}((m+1)K_{X}))).
\]
Then we define 
\[
K_{m+1} = \sum_{i=0}^{N_{m+1}}\mid\!\sigma_{i}^{(m+1)}\!\!\mid^{2}
\]
and 
\[
h_{m+1}:= 1/K_{m+1}
\]
inductively.  

\begin{theorem}\label{main2}
 Let $X$ be a smooth projective $n$-fold with ample canonical 
bundle.  Let $A$ and $\{ h_{m}\}_{m \geqq 1}$ be as above.
Then 
\[
h_{\infty} := \liminf_{m\rightarrow\infty} \sqrt[m]{(m!)^{n}\cdot h_{m}}
\]
is a hermitian metric on $K_{X}$ such that 
\[
\omega_{\infty} := \sqrt{-1}\Theta_{h_{\infty}}
\]
is a K\"{a}hler-Einstein current $\omega_{E}$ on $X$ as in Theorem \ref{KE}. 
$\square$
\end{theorem}
The proof of Theorem \ref{main2} is essentially the same as the one of 
Theorem \ref{main}.  

Similarly to Lemma \ref{lower}, we obtain the following lower estimate. 
\begin{lemma}\label{lower2}
\[
\limsup_{m\rightarrow\infty}\frac{1}{(m!)^{n/m}}\int_{X}(K_{m})^{\frac{1}{m}}
\geqq (2\pi)^{-n}dV_{E}
\] 
holds. $\square$
\end{lemma}
Lemma \ref{lower2} can be obtained just as in the proof of  Lemma \ref{lower}.  In the proof of Lemma \ref{lower}, we have considered 
all the points on $X$, but here we only need to consider points on 
\[
X^{\circ}:= \{x \in X \mid \mbox{$\Phi_{\mid mK_{X}\mid}$ is an embedding 
on a neighbourhood of $x$ for some positive integer $m$}\}.  
\]
$X^{\circ}$ is the locus where $\omega_{E}$ is $C^{\infty}$ strictly 
positive form.
In fact the proof of Lemma \ref{lower} is essentially local (as $m$ tends to infinity). \vspace{3mm} \\ 
For the upper estimate, we set 
\[
\mu (X,K_{X}) = n!\cdot\limsup_{m\rightarrow\infty}m^{-n}\dim H^{0}(X,{\cal O}_{X}(mK_{X})).
\]
Then by the same manner as the proof of  Lemma \ref{upper}, we obtain the 
following lemma. 
\begin{lemma}\label{upper2}
\[
\limsup_{m\rightarrow\infty}\frac{1}{(m!)^{n/m}}(K_{m})^{\frac{1}{m}} \leqq \frac{1}{n!}\mu (X,K_{X}). 
\]
holds. $\square$
\end{lemma}
Let $\pi_{m} : X_{m}\longrightarrow X$ be the resolution of $\mbox{Bs}\mid m!K_{X}\mid$.  
Let 
\[
\pi_{m}^{*}\mid\! m!K_{X}\!\mid = \mid\! P_{m}\!\mid + F_{m}
\]
be the decomposition of $\pi_{m}^{*}\mid\!\! m!K_{X}\!\!\mid$ into the 
free part $\mid\!\!P_{m}\!\!\mid$ and the fixed component $F_{m}$. 
By Fujita's theorem 
(\cite[p.1, Theorem]{f}), we see that 
\[
\lim_{m\rightarrow\infty}\frac{P_{m}^{n}}{(m!)^{n}} 
= \mu (X,K_{X})
\]
holds.  
Then by the construction of $\omega_{E}$ (cf. Section \ref{exist}), we see that 
\[
\frac{1}{(2\pi)^{n}}\int_{X}dV_{E} = \frac{1}{n!}\mu (X,K_{X})
\]
holds.  Hence combining Lemmas \ref{lower2} and \ref{upper2}, 
we obtain that 
\[
\limsup_{m\rightarrow\infty}\frac{1}{(m!)^{n/m}}\int_{X}(K_{m})^{\frac{1}{m}}
= (2\pi)^{-n}\int_{X}dV_{E}
\]
and 
\[
\limsup_{m\rightarrow\infty}\frac{1}{(m!)^{n/m}}(K_{m})^{\frac{1}{m}}
= (2\pi)^{-n}dV_{E}
\]
hold.  This completes the proof of Theorem \ref{main2}. $\square$

\begin{corollary}
Let $X$ be a smooth projective variety of general type.
Then the maximal K\"{a}hler-Einstein current (cf. Definition \ref{maximal})
is unique.  $\square$
\end{corollary}
\section{Proof of Theorems \ref{positivity} and \ref{nakano2}}

In this section we shall prove  Theorems \ref{positivity},\ref{nakano2} 
by using Theorems \ref{main},\ref{main2}.   
Roughly speaking, Theorems \ref{main},\ref{main2} imply that what we can say about Bergman kernels also holds for K\"{a}hler-Einstein volume forms. 
\vspace{5mm} \\

\noindent{\em Proof of Theorems \ref{positivity},\ref{nakano2}}. Let $A$ be a sufficiently ample line bundle on $X$ 
and let $h_{0}$ be a $C^{\infty}$ hermitian metric with strictly 
positive curvature.

Then for every $s\in S^{\circ}$, we define the dynamical system
of the Bergman kernels $\{ K_{m,s}\}$ on the fiber $X_{s}:= f^{-1}(s)$
as in Section \ref{m2}. 
Then we see that the hermitian metric
\[
h_{m}\mid X_{s} = 1/K_{m,s}
\]
on $A + mK_{X/S}\mid X^{\circ}$ has semipositive curvature  by Theorem \ref{v}.
And it extends to a singular hermitian metric on $A + mK_{X/S}$ with semipositive curvature as in Theorem \ref{v}.
Then by Theorem \ref{main2}, we see that $h_{E}$ is a singular hermitian metric
on $K_{X/S}$ with semipositive curvature current. 
This completes the proof of Theorem \ref{positivity}. 
Then by Theorem \ref{nakano}, we complete the proof of Theorem \ref{nakano2}. 
$\square$ 
\begin{remark} In \cite{tu5}, I have proved the Nakano semipositivity of 
$f_{*}{\cal O}_{X}(mK_{X/S})$ similar to Theorem \ref{nakano2} even when a general fiber is of non-general type. But in this case the metric does not come 
from  K\"{a}hler-Einstein metrics. 
$\square$
\end{remark}

Author's address\\
Hajime Tsuji\\
Department of Mathematics\\
Sophia University\\
7-1 Kioicho, Chiyoda-ku 102-8554\\
Japan \\
e-mail address: tsuji@mm.sophia.ac.jp


\begin{thebibliography}{99}
\bibitem[A]{a} Aubin, T.: Equation du type Monge-Amp\`{e}re sur les variet\'{e}
 k\"{a}hlerienne compactes, C.R. Acad. Paris {\bf 283} (1976), 459-464.
\bibitem[B-T]{b-t} Bedford, E. and  Taylor, B.A., A new capacity of plurisubharmonic  functions, Acta Math. {\bf 149} (1982), 1-40.

\bibitem[B1]{b1} Berndtsson, B.: Subharmonicity properties of the Bergman kernel and some other functions associated to pseudoconvex domains, math.CV/0505469 (2005). 
\bibitem[B2]{b2} Berndtsson, B.: Curvature of vector bundles and subharmonicity
of vector bundles, math.CV/050570 (2005).
\bibitem[B3]{b3} Berndtsson, B.:  Curvature of vector bundles associated to holomorphic fibrations, math.CV/0511225 (2005). 
\bibitem[C]{c} Catlin, D.: The Bergman kernel and a theorem of Tian, 
Analysis and geometry in several complex variables (Katata 1997),  
1-23, Trends in Math.,
Birkh\"{a}user Boston, Boston MA. (1999). 
\bibitem[D-P-S]{d-p-s}  Demailly, J.P.-Peternell, T.-Schneider, M. : 
Pseudo-effective line bundles on compact K\"{a}hler manifolds, 
math. AG/0006025 (2000).
\bibitem[Do]{do} Donaldson, S.K.: Scalar curvarure and projective embeddings I,
Journal of Differential Geom. {\bf 59} (2001), 479-522.
\bibitem[F]{f} Fujita, T.: Approximating Zariski deecomposition
of big line bundle , Kodai Math. J. {\bf 17} (1994), 1-4. 
\bibitem[Ka]{ka1} Kawamata, Y.: Kodaira dimension of Algebraic giber spaces over curves, Invent. Math. {\bf 66} (1982), pp. 57-71.
\bibitem[Kob]{kob} Kobayashi, R.: Existence of K\"{a}hler-Einstein metrics
on an open algebraic manifold, Osaka J. of Math. {\bf 21} (1984), 399-418.
\bibitem[Kr]{kr} Krantz, S.: Function theory of several complex variables, 
John Wiley and Sons (1982).  
\bibitem[M-Y]{m-y} Maitani, and Yamaguchi, S.: 
Variation of Bergman metrics on Riemann surfaces, Math. Ann. {\bf 330} (2004) 477-489.
\bibitem[N]{n}Nadel, A.M.: Multiplier ideal sheaves and existence of K\"{a}hler-Einstein
metrics of positive scalar curvature, Ann. of Math. {\bf 132}(1990),549-596.
\bibitem[O-T]{o-t}Ohsawa, T and  Takegoshi K.: $L^{2}$-extension of holomorphic
functions, Math. Z. 195 (1987),197-204.
\bibitem[O]{o} Ohsawa, T.: On the extension of $L^{2}$ holomorphic functions V,
effects of generalization, Nagoya Math. J. {\bf 161}(2001) 1-21.  
\bibitem[Ti]{ti}Tian, G.: On a set of polarized K\"{a}hler metrics on algebraic
manifolds, Jour. Diff. Geom. {\bf 32}(1990),99-130. 
\bibitem[Tr]{tr} Trudinger, N.S.: Fully nonlinear elliptic equation under 
natural structure conditions, Trans. A.M.S. {\bf 272} (1983), 751-769.
\bibitem[S]{s} Sugiyama, K.: Einstein-Kahler metrics on minimal varieties of general type and an inequality between Chern numbers. 
Recent topics in differential and analytic geometry, 417--433, 
Adv. Stud. Pure Math., {\bf 18}-{\rm I}, 
Academic Press, Boston, MA, 1990. 
\bibitem[T0]{tu0} Tsuji H.: Existence and degeneration of Kahler-Einstein metrics on minimal algebraic varieties of general type. Math. Ann. {\bf 281} (1988), no. 1, 123--133. 
\bibitem[T1]{tu}Tsuji H.: Analytic Zariski decomposition, Proc. of Japan Acad.
{\bf 61}(1992), 161-163.
\bibitem[T2]{tu2} Tsuji, H.:  Existence and Applications of Analytic Zariski Decompositions, Trends in Math., Analysis and Geometry in Several Complex Variables(Katata 1997), Birkh\"{a}user Boston, Boston MA.(1999), 253-272.
\bibitem[T3]{tu3}Tsuji, H.: Deformation invariance of plurigenera, Nagoya Math. J. {\bf 166} (2002), 117-134.
\bibitem[T4]{tu4} Tsuji, H.:  Refined semipositivity and Moduli of canonical models, preprint (2005). 
\bibitem[T5]{tu5} Tsuji, H: Variation of Bergman kernels of adjoint line bundles,math.CV/0511342 (2005).
\bibitem[Y1]{y} Yau, S.-T.: On the Ricci curvature of a compact K\"{a}hler manifold and the complex Mong\'{e}-Amp\'{e}re equation,  Comm. Pure  Appl. Math. {\bf 31} (1978),339-441.
\bibitem[Y2]{y2} Yau, S.-T.: A general Schwarz lemma for K\"{a}hler manifolds,
Amer. J. of Math. {\bf 100} (1978), 197-203. 
\bibitem[Ze]{ze} Zelditch, S.: Sz\"{o}go kernel and a theorem of Tian, 
International Reserch Notice 6 (1998),
317-331.
\bibitem[Z]{z} Zhang, S.: Heights and reductions of semistable varieties, 
Compositio Math. {\bf 104}(1996), 77-105. 
\end{thebibliography}
\end{document}